\documentclass[12pt]{article}
\usepackage[english]{babel}
\usepackage[applemac]{inputenc}
\usepackage{amsmath,amssymb,mathrsfs}
\usepackage{graphicx,pifont,hyperref}
\usepackage{amsthm,amscd}
\usepackage[left=25mm, right=25mm, top=30mm, bottom=30mm]{geometry}

\newtheorem{theo}{Theorem}[section]
\newtheorem{prop}[theo]{Proposition}

\newtheorem{defi}[theo]{Definition}
\newtheorem{coro}[theo]{Corollary}
\newtheorem{rem}[theo]{Remark}
\newtheorem{prob}[theo]{Problem}
\newcommand{\eps}{\varepsilon}
\newcommand{\set}[1]{\left\{#1\right\}}
\newcommand{\bbR}{\mathbb{R}} 
\newcommand{\ie}{\emph{i.e.}\ }

\newcommand{\norme}[1]{\parallel #1\parallel}
\newcommand{\Alg}{\mathcal{A}}
\newcommand{\cste}{\mathrm{c}^{\text{ste}}}
\newcommand{\marginal}[1]{\marginpar{\raggedright\tiny #1}}
\newcommand{\PC}[1]{\marginal{PC:\;#1}}

\newcommand{\LG}[1]{\marginal{LG:\;#1}}

\newcounter{NumRem}

\newenvironment{algorithme}[1][]
{
\vspace{3mm}

\noindent{\bf Algorithm.}\ttfamily \\
\noindent}
{
\vspace{3mm}\\
}

\newenvironment{intuition}[1][]
{\begin{center}
\begin{minipage}{110mm}
\it
}
{
\end{minipage}
\end{center}
}

\title{Efficient estimation of the cardinality of large data sets}
\author{Philippe Chassaing and Lucas Gerin}
\begin{document}
\maketitle
\abstract{The estimation of the number  of distinct elements in a sequence of words, under strong constraints 
coming from applications to database analysis and network routing, is the problem we address in this article. 
Giroire has recently proposed a solution, under the form of a probabilistic algorithm 
that uses statistical properties of uniform random variables in [0,1]. 
Our objective here is twofold :

First, the analysis of this algorithm within 
the framework of Kullback information and estimation theory allows us to reinterpret
a lower bound due to Indyk \& Woodruff as a consequence of well-known inequalities.

Second, we show that a slight modification of the Giroire algorithm returns an estimation
whose accuracy is optimal, among algorithms based on order statistics.}

\noindent\emph{NB: This paper is the extended version of \cite{CGEstimation}}
\section{Introduction}
\newcommand{\moins}{\smallsetminus}

Let $\mathbf{y}_n=(y_1,y_2,\dots,y_n)$ be a sequence of elements of a finite set $\mathcal{C}$.
We wish to design an algorithm that computes $\theta=\mathrm{card}\set{y_1,y_2,\dots,y_n}$.
This question is sometimes referred to as the \emph{distinct elements problem}.
Consider the following naive algorithm:
\begin{algorithme}
Initialize Dictionary=$\{\}$\\
For j=1 to n\\
$\phantom{12}$  Look for $y_j$ in Dictionary \\
$\phantom{12}$  If $y_j$ is in Dictionary, do nothing\\
$\phantom{12}$  otherwise, add $y_j$ to Dictionary\\
next j\\
Return the size of Dictionary 
\end{algorithme}
During the execution of this algorithm, each word must be stored on the disk, so that the  memory requirement cannot be less than linear in $\theta$.
For each $y_j$, a query in \verb|Dictionary| is needed, that costs at least  
$\mathcal{O}(\log \theta)$ elementary operations. For a number of applications, these linear-space and log-time complexity lower bounds are not satisfactory, but, for the time being, they cannot be improved (see \cite{AMS} for a discussion). For instance, finding the number of distinct elements in a given sequence is a main issue in network routing, where huge  data sets have to be handled : according to  \cite{Gi}, typically, at a given node of a network, packets arrive every 60 nanoseconds,  approximately, but hardware limitations make impossible to process more than 100 elementary operations on each packet, and the size of $\mathcal{C}$ does not allow to store every word.
\begin{prob}\label{Prob}
Does there exist an algorithm $\Alg$ that returns the number of different 
elements in $\mathbf{y}_n$, while statisfying the two following constraints:
\begin{enumerate}
\item $\Alg$ uses at most $M$ bits of memory  ;
\item each data $y_i$ is treated in one pass: a few operations are processed (possibly modifying the $M$ bits), then $y_i$ is irreversibly erased.
\end{enumerate}
\end{prob}
These constraints are too strong to allow an exact solution of the distinct elements problem, for a memory of $M$ bits counts  at most $2^M$ distinct elements. Partly following \cite{Mo}, Flajolet and Martin \cite{FM} proposed a probabilistic algorithm that overcomes this limitation by returning only an approximation of $\theta$. 
\subsection{Approximate counting}
We randomize  Problem \ref{Prob} through the use of \emph{hash functions} :
\begin{defi}
Given a \emph{typical} sequence of distinct words, a \emph{hash function} $h:\ \mathcal{C}\to [0,1]$
returns a sequence of random numbers, i.e.  a sequence of numbers that behaves as the realization of a sequence of independent random variables, uniform
on $[0,1]$.
\end{defi}
\noindent We assume from now on
that we are given this idealized version $h$ of a hash function, and that the \textit{set}
\begin{displaymath}
\mathbf{X}=\{X_1,\dots,X_n\},
\end{displaymath}
where $X_i=h(y_i)$,
is distributed like a \textit{set} of $\theta$ independent realisations of a uniform random variable
on $[0,1]$. The design of \emph{good} hash functions is discussed for example in
Knuth's book \cite{Knu3}.

The algorithm proposed by Flajolet and Martin \cite{FM}  is based on the distribution of 
some patterns of the dyadic representation of the $X_i$'s.
Their work has revealed the following phenomenon.
\begin{intuition}
It is possible to recover (an approximate value of) $\theta$,
using only a (small and) constant memory $M$.
\end{intuition}
Expectedly, the accuracy of the approximation  of $\theta$  increases with $M$.
Indyk and Woodruff provide  the following theoretical bound for the accuracy reached with the help of a $M$-bits memory.
\begin{prop}[Indyk-Woodruff \cite{IW}]\label{Th:IW}
For fixed $\eps,\delta>0$, one says that a probabilistic algorithm $\Alg(\eps,\delta)$-approximates $\theta$
if it returns a value $\hat{\theta}$ such that
$$
\mathbb{P}(|\hat{\theta}-\theta|>\theta\eps)\leq\delta.
$$
Let $\Alg$ a one-pass algorithm that $(\eps,\delta)$-approximates $\theta$, 
and assume that $\eps=\mathcal{O}(a^{-\frac{1}{9}})$,
in which $a$ is the number of elements in $\mathcal{C}$.
Let $M$ be the memory required by $\Alg$ ($M$ is expressed in bits). Then
$$
\eps^{-1}= \mathcal{O}\left(\sqrt{M}\right).
$$
\end{prop}
The proof of this bound is derived from learning theory (mainly, from bounds given  in terms of VC-dimensions of well-choosen sets \cite{Vap}), and requires a fine study of the geometry
of some $\ell^1$ or $\ell^2$ spaces.

Assume that $\Alg$ returns an unbiased estimation of $\theta$ (\ie $\mathbb{E}[\hat{\theta}]=\theta$). From the Bienaim-Tchebychev inequality and Proposition \ref{Th:IW}, we obtain
\begin{align*}
\mathbb{P}(|\hat{\theta}-\theta|>\theta\eps)
&\leq \frac{\mathrm{Var}(\hat{\theta})}{(\theta\eps)^2}\\
&\leq \cste\ \frac{\mathrm{Var}(\hat{\theta})}{\theta^2}\ M,
\end{align*}
leading to the following lower bound :
$$
\mathrm{Var}(\hat{\theta})\geq \cste\ \frac{\theta^2\,\mathbb{P}(|\hat{\theta}-\theta|>\theta\eps)}{M}.
$$
Thus, an interpretation of Proposition \ref{Th:IW} is that the variance of $\hat{\theta}$ cannot decrease faster than $\frac{\theta^2\,\mathbb{P}(|\hat{\theta}-\theta|>\theta\eps)}{M}$. The main goal of the present article is to analyse the algorithm  {\scshape MinCount}, proposed by Giroire \cite{Gi}, within the framework of estimation theory. For this algorithm, and other algorithms based on \emph{order statistics}, we shall show that a more precise lower bound of $\theta^2/M$ appears. It does not follow from geometric considerations but as a consequence of well-known inequalities in estimation theory. We show furthermore that a slightly improved variant of {\scshape MinCount} achieves this bound.
\subsection{The {\scshape MinCount} algorithm}
First, we describe a simplified version of {\scshape MinCount} \cite{Gi}, with parameter
 $k$ ($k$ is a given  integer, not smaller than 2). Each word $y_i$ is hashed, the corresponding value
is compared to the $k$ smallest values already observed. As usual,
$X_{(1)}=\min_{j\leq n} X_j$ stands for  the smallest $X_j$, and $X_{(k)}$ for the
$k$-th smallest.
\begin{algorithme}
Set MIN[1]=MIN[2]=...=MIN[k]=1\\
For j=1 to n,\\
$\phantom{123}$ Compare $X_j:=h(y_j)$ with MIN[1:k]\\
$\phantom{123}$ Update the vector (MIN[1]$\leq$MIN[2]$\leq$ ...$\leq$ MIN[k]) of the $k$ smallest \\
$\phantom{123}$ values of the sequence $(X_k)_{1\le k\le j-1}$ according to the value  of $X_j$ \\
next j\\
Return $\frac{k-1}{\mbox{{\ttfamily MIN[k]}}}$.
\end{algorithme}

This simplified algorithm satisfies the constraints of Problem \ref{Prob}, provided that  $k=\mathcal{O}\left(M\right)$ (we are more precise below). Indeed, {\scshape MinCount} processes each word $y_{i}$ using a single pass,
and throughout the execution, only $k$ real numbers are kept in memory.
To see why $(k-1)/X_{(k)}$ gives a good estimation of $\theta$, recall from \cite[pages 8-13]{Dav} the density of probability of the $k$-th order statistic:
\begin{displaymath}
\mathbb{P}(X_{(k)}\in[t,t+dt))=\theta\binom{\theta-1}{k-1}t^{k-1}(1-t)^{\theta-k}dt.
\end{displaymath}
It follows that $\mathbb{E}[1/X_{(1)}]=+\infty$, but as soon as $k\geq 2$,
\begin{equation}\label{Eq:EstNaif}
\mathbb{E}[1/X_{(k)}]=\theta\binom{\theta-1}{k-1}\mathrm{B}(k-1,\theta-k+1)
=\frac{\theta}{k-1},
\end{equation}
where $\mathrm{B}$ is the Euler beta function.

If a number of the unit interval is stored with a precision of $2^{-s}$, its storage requires $s$ bits, and the available number of bits,  $M$,  allows to store $k=M/s$ numbers. 
In \cite{Gi},  {\scshape MinCount}  is tuned  in two directions : 
\begin{enumerate}
\item the interval $[0,1]$ is split into  $m$ sub-intervals $[0/m,1/m)$, $[1/m,2/m)$,  $\dots$ ; the algorithm returns the $k$-th smallest value lying in each of these intervals ; 
\item rather than $x\mapsto (k-1)/x$, three different functions are proposed (depending on the parameters $k,m$).
\end{enumerate}
The division of $[0,1]$ in $m$ intervals, called \emph{stochastic averaging} in \cite{FM}, allows to obtain a sample of $m$ copies of $X_{(k)}$ almost at the same cost as one copy : when $m=1$, each $X_j$ has to be compared with the $k=M/s$ smallest values stored at time $j$, while, when $m\neq 1$, one has first to find $i$ such that 
$$\frac{i-1}{m}\leq X_j< \frac{i}{m},$$
which can be done at almost no cost\footnote{It is a simple truncature if $m$ is a power of $2$.}, then $X_i$ is compared to the $\tilde{k}=M/(m\tilde{s})$ smallest values lying in $[(j-1)/m,j/m)$, stored at time $j$.

Given the number $M$ of bits, we see that a trade-off has to be made between $m$ large (fewer comparisons, and a larger sample, but $\tilde{k}$ or $\tilde{s}$ smaller) and $m$ small ($k$ and $s$ are larger). In what follows, we shall study the impact of $k$ and $m$ on the accuracy of the estimation of $\theta$, and we shall consider that $$M=k\times m,$$
thus not taking the impact of $s$ into account. Here is the algorithm {\scshape MinCount}, as given in \cite{Gi}.
\begin{algorithme}
Fix two integer parameters $k\geq 2,m\geq 1$.\\
Set $Z_{(p),i}=\frac{i}{m}$ for each $i\leq m,p\leq k$.\\
For j=1 to N\\
$\phantom{1}$ $Z_j=h(y_j)$.\\
$\phantom{1}$ Let $i$ such that $Z_j\in [\frac{i-1}{m},\frac{i}{m})$.\\
$\phantom{1}$ Update the vector $(Z_{(1),i},\dots,Z_{(k),i})$ of the $k$ smallest values in $[\frac{i-1}{m},\frac{i}{m})$.\\
next $j$.\\
For each $p,i$, set $X_{(p),i}=Z_{(p),i}-\frac{i-1}{m}$.\\
\\
Return a function $\hat{\xi}=\hat{\xi}(X_{(l),i};1\leq i\leq m;1\leq l\leq k)$.
\end{algorithme}
This algorithm fulfills the constraints of Problem \ref{Prob}, and
the memory needed, when expressed in bits, is linear in $M:=k\times m$.

The distinct elements problem is now reduced to a statistical problem : given a $k\times m$-\emph{sample}
$$
\Xi_{k,m}=\left(X_{(1),1},\dots,X_{(k),1},\dots,X_{(1),m},\dots,X_{(k),m}\right)
$$
the distribution of which depends on an unknown parameter $\theta$, one has to find the best \emph{estimation} of $\theta$.
\begin{figure}
\begin{center}
\includegraphics[width=80mm]{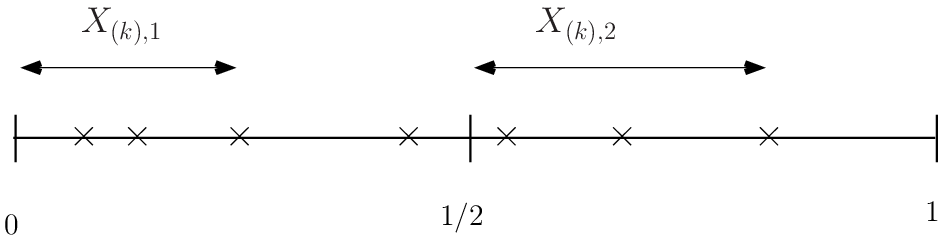}
\end{center}
\caption{An example of the sample $(X_{(k),1},\ldots,X_{(k),m})$,
 with $\theta=7,m=2,k=3$. The crosses represent the 7 hashed values.}
\label{Fig:Buck}
\end{figure}
Giroire \cite{Gi} proposes three different estimators $\xi_1,\xi_2,\xi_3$, each depending on $\Xi_{k,m}$ :\PC{y'a un bleme avec la taille des captions}
\begin{eqnarray*}
\xi_1 &=& \frac{k-1}m\ \sum_{i=1}^m\frac{1}{X_{(k),i}},
\\
\xi_2 &=& C(k,m)\ \left(\sum_{i=1}^m\frac{1}{\sqrt{X_{(k),i}}}\right)^2,
\\
\xi_3 &=& \biggl(\frac{\Gamma(k-1/m)}{\Gamma(k)}\biggr)^{-m}\exp\left(-\frac{1}{m} \sum_{i=1}^m\log X_{(k),i}\right).
\end{eqnarray*}
According to \cite{Gi}, the $\xi_i$'s are asymptotically  unbiased: $\mathbb{E}[\xi_i]\sim \theta$ as  $\theta\to+\infty$. To our knowledge, for the distinct elements problem, Giroire's algorithms leading to the $\xi_i$ were the best one-pass  algorithms, i.e. producing the estimation with the lowest quadratic error, with $\xi_3$ ahead of the 2 others. One of the goals of this paper is to show that the best (with respect to quadratic error) estimation of $\theta$ based on $\Xi_{k,m}$, is given by
\begin{displaymath}
\hat{\xi}=\frac{km-1}{\sum_{i=1}^{m} X_{(k),i}}.
\end{displaymath}
\begin{rem}
Given $X_{(k),i}=x<\tfrac1m$, the random variables $(X_{(1),i},\dots,X_{(k-1),i})$ are uniformly distributed on $[0,x]$, i.e. their conditional distribution does not depend on $\theta$. In other words, given $X_{(k),i}$, the knowledge of the $k-1$ observations $(X_{(1),i},\dots,X_{(k-1),i})$ does not  bring any additional information on $\theta$, so it comes as no surprise that $\hat{\xi}$, or even  $\xi_1,\xi_2,\xi_3$, depend only on $X_{(k),i}$.
\end{rem}

\noindent{\bf Structure of the paper.} The number of values falling in each of the $m$ subintervals follows a multinomial distribution, and, conditionally, given the number of values falling in each of the $m$ subintervals, the distribution of the  $X_{(k),i}$'s is the distribution of the $k$-th order statistic of a uniform random sample with random (binomial) size. The reader surely admits that this description does not sound very tractable. In the next section, we shall thus study the asymptotic distribution of the $X_{(k),i}$'s for $\theta$ large, for this asymptotic distribution is much simpler than the true distribution. Then we shall prove that $\hat{\xi}$ is the best estimator, provided that the sample $\Xi_{k,m}$ follows the asymptotic distribution, and we shall explain the lower bound given by \cite{IW}. In the last Section, we discuss the optimality of $\hat{\xi}$ when $\Xi_{k,m}$ follows its actual distribution.
\section{The best estimation in the asymptotic model}
\subsection{The asymptotic model}
\label{sec:asympt}
First we describe the asymptotic behavior of the $X_{(k),i}$'s. Recall that a random variable $\Gamma_{k,\theta}$ is said to follow the Gamma distribution with  parameters $k$ and $\theta$ if
\begin{displaymath}
\mathbb{P}(\Gamma_{k,\theta}\in[t,t+dt))
=\frac{t^{k-1}}{\Gamma(k)}\theta^{k}e^{-\theta t}\mathbf{1}_{t\geq 0}\,dt.
\end{displaymath}
Furthermore, the vector $(m_{\theta,i})_{1\le i\le k}$ of the $k$ smallest among $\theta$ i.i.d. random variables, uniform on $[0,1]$, satisfies
\begin{displaymath} 
\theta\, (m_{\theta,i})_{1\le i\le k} \xrightarrow[\theta\to\infty]{\text{(law)}}T_{k}Y,
\end{displaymath}
in which the $k$ components of the column vector $Y$ are i.i.d. and follow the exponential distribution with expectation $1$, i.e. the Gamma distribution with  parameters $1$ and $1$, and $T_{k}$ is the $k\times k$ matrix with ones on and below the diagonal, and zeroes above the diagonal. As a consequence, 
\begin{displaymath} 
\theta\, m_{\theta,k} \xrightarrow[\theta\to\infty]{\text{(law)}}\Gamma_{k,1}.
\end{displaymath}
Since there are approximately $\theta/m$ elements in each subinterval, and they are distributed as i.i.d. random variables, uniform on $[0,1/m]$, we can expect that 
\begin{equation}\label{Eq:CvGamma}
(\theta X_{(k),1},\ldots,\theta X_{(k),m})
\xrightarrow[\theta\to\infty]{\text{(law)}}(\Gamma^{(1)},\ldots,\Gamma^{(m)}),
\end{equation}
in which the $\Gamma^{(i)}$'s are $m$ i.i.d. r.v. with common distribution Gamma $(k,1)$.  Also, we can expect $\theta\,\Xi_{k,m}$  to be distributed, asymptotically, as $m$ independent copies of $T_{k}Y$.
Since
\begin{displaymath}
\frac{1}{\theta}\Gamma_{k,1}\stackrel{\text{(law)}}{=}\Gamma_{k,\theta},
\end{displaymath}
equation \eqref{Eq:CvGamma} roughly says that, when $\theta$ goes large, the 
 $X_{(k),m}$ behave as $m$ independent random variables with distribution Gamma$(k,\theta)$.
We shall not prove \eqref{Eq:CvGamma}, for we need only a more specific result : the quadratic error is asymptotically the same when the $X_{(k),i}$'s are replaced by $m$
independent Gamma random variables. As we shall see later,  
\[\mathbb{E}[(\hat{\xi}-\theta)^2]=\mathcal{O}(\theta^2).\]
Thus we only need to compare $\hat{\xi}$ to  functions $f(\Xi_{k,m})$ such that $\mathbb{E}[(f(\Xi_{k,m})-\theta)^2]$ is not too large.
\begin{prop}\label{Prop:Conv}
Let $f$ be a continuous function: $\mathbb{R}_+^m \backslash \{0\}\to\mathbb{R}_+$.
We assume that outside some neighbourhood of $0$,  $f$ is bounded, while, in the neighbourhood of $0$,  there exists $r>0$ such that
\begin{equation}\label{Eq:zeropole}
|f(x)| =\mathcal{O}\left(\frac{1}{\parallel x\parallel^r}\right).
\end{equation}
Then, for any $\eps>0$, there exists $c_\eps>0$ such that, for $\theta$ large enough,
\LG{Terme exponentiellement petit en $\theta^2$ pour tre tranquille}
\begin{multline}\label{Eq:Cv:Prop}
\left|\mathbb{E}\left[\left(f(X_{(k),1},\ldots,X_{(k),m})-\theta\right)^2\right]
-\mathbb{E}\left[\left(f(\Gamma^{(1)},\dots,\Gamma^{(m)})-\theta\right)^2\right]\right|\\
\leq c_\eps \theta^{\eps -1}\mathbb{E}\left[\left(f(X_{(k),1},\ldots,X_{(k),m})-\theta\right)^2\right]
+4m\theta^2 \exp(-\theta/(2m)^2).
\end{multline}
in which the $\Gamma^{(i)}$'s are i.i.d. random variables with law Gamma$(k,\theta)$.
\end{prop}
Consequently, we shall first assume that $\theta$ is large enough, and we shall replace the  $\set{X_{(k),i},i=1,\dots,m}$ with their limits $\set{\Gamma^{(i)},i=1,\dots,m}$.  We know that the $X_{(k),i}$'s are large only if $\theta$ is small, thus the assumption of boundedness of $f$ outside a neighbourhood of 0 is not really binding for an estimator of $\theta$. Similarly, a good estimation has to be moderately large when    the $X_{(k),i}$'s are small, thus a good estimation fulfills necessarily \eqref{Eq:zeropole}.   The proof of Proposition \ref{Prop:Conv} is postponed to Section \ref{Sec:Preuve}. 
\subsection{Lehmann-Scheff and Cramr-Rao inequalities}
We are now given a \emph{sample} $(\Gamma^{(1)},\dots,\Gamma^{(m)})$ of $m$ independent
Gamma r.v. with parameters $(k,\theta)$. We assume that $k$ is known, and we want to estimate
the unknown parameter $\theta$.
We proceed by maximum likehood estimation: let $f_\theta : \mathbb{R}_+^m\to\mathbb{R}_+$
be the density of the $m$-sample :
\begin{align*}
f_\theta(x_1,\dots,x_m)&=\ \theta^{km}\ \exp(-\theta \textstyle{\sum_i x_i})\ \prod_i \frac{{x_i}^{k-1}}{\Gamma(k)}.
\end{align*}
We are to compute
\begin{displaymath}
\hat{\theta}:=\mathrm{argmax}_{\theta >0}\ \ln\left(f_\theta(\Gamma^{(1)},\dots,\Gamma^{(m)})\right),
\end{displaymath}
thus we have to solve :
\begin{eqnarray*}
0
&=&
\frac{\partial}{\partial \theta}\ \ln\left(f_\theta(\Gamma^{(1)},\dots,\Gamma^{(m)})\right)
\\
&=&
\frac{km}{\theta}\ -\ \sum_i \Gamma^{(i)},
\end{eqnarray*}
leading to $\tilde{\theta}=\frac{km}{\sum_i \Gamma^{(i)}}$.
It turns out that this estimator is biased: 
\[\mathbb{E}[\tilde{\theta}]=\theta \frac{km}{km-1}.\]
The next Proposition fixes the problem :
\begin{prop}
\label{calculdebase}
Set
$$\hat{\xi}(x_1,\dots,x_m)=\frac{km-1}{x_1+\dots +x_m}.$$
For any $\theta$,
\begin{align*}
\mathbb{E}[\hat{\xi}(\Gamma^{(1)},\dots,\Gamma^{(m)})] & =\theta,\\ 
\text{Var}(\hat{\xi}(\Gamma^{(1)},\dots,\Gamma^{(m)})) & =\frac{\theta^2}{km-2}.
\end{align*}
\end{prop}
\begin{proof}
By stability of the class of Gamma distributions under convolution, 
\begin{displaymath}
\gamma:=\Gamma^{(1)}+\dots +\Gamma^{(m)} 
\end{displaymath}
is Gamma distributed with parameters $km$ and $\theta$. Thus,
\begin{align*}
\mathbb{E}\left[\frac{km-1}{\gamma}\right]
&=
\ \frac{(km-1)\theta^{km}}{\Gamma(km)}\  \int_0^{+\infty}\ {x}^{km-2} e^{-\theta x} dx
\\
&=
\ \frac{(km-1)\Gamma(km-1)}{\Gamma(km)}\  \theta\ =\ \theta.
\end{align*}
The variance is obtained through  similar computations.
\end{proof}
Combined with Proposition \ref{Prop:Conv}, Proposition \ref{calculdebase} entails that
\begin{coro}
\label{quaderror}
$$
\mathbb{E}\left[\left(\hat{\xi}(X_{(k),1},\ldots,X_{(k),m})-\theta\right)^2\right]=
\frac{\theta^2}{km-2}+\mathrm{o}(\theta^2).
$$
\end{coro}
Corollary \ref{quaderror}  calls for two remarks:
\begin{enumerate}
\item the asymptotic variance $\theta^2/(km-2)$ is indeed smaller than the variance of the estimators proposed in \cite{Gi}, and sets a new record for the lower quadratic error for a one-pass algorithm requiring bounded memory ;
\item since the memory required is linear in $km$, the quadratic error obtained in Corollary \ref{quaderror}  is consistent with the theoretical result given in \cite{IW}.
\end{enumerate}

Let us now recall a few definitions from estimation theory.
Given a $m$-sample $(X_1,\ldots,X_m)$, whose law is denoted $P_\theta$, any random variable  $S=S(X_1,\ldots,X_m)$ is called a \emph{statistic}. Here, $(X_1,\ldots,X_m)=(\Gamma^{(1)},\dots,\Gamma^{(m)})$ and we shall focus on the statistic $S=\sum_i \Gamma^{(i)}$. First, we check that $S$ fulfills two conditions with deep connections with accuracy : \emph{sufficiency} and \emph{completeness}.
\begin{defi}\label{Def:Ex}
Assume that $P_\theta$ admits a density $f_{\theta}$ with respect to the Lebesgue measure :
\begin{enumerate}
\item a statistic $S$ is said to be \emph{sufficient} for $\theta$ if $f_\theta$ can be written
\begin{displaymath}
f_\theta(x_1,\dots,x_m)\ =\ g(S(x_1,\dots,x_m),\theta)\ h(x_1,\dots,x_m),
\end{displaymath}
in which $g,h$ are two non-negative measurable functions, $h$ not depending on $\theta$.
\item a statistic $S$ is said to be \emph{complete} if, for any measurable function $\phi$,
$$
\left\{\forall \theta,\ \mathbb{E}[\phi(S)]=0\right\} \Rightarrow
\left\{\forall \theta,\ \{\phi(S)\equiv 0, P_\theta-\mbox{p.s.}\}\right\}.
$$
\end{enumerate}
\end{defi}
The next criterion ensures sufficiency and 
completeness : 
\begin{prop}[see \cite{Lind},Th.16 Chap.7]\label{Prop:Critere}
Assume that $f_\theta$ can be written 
\begin{displaymath}
f_\theta(x)=h(\mathbf{x})B(\theta)\exp\left(Q(\theta)R(x)\right),
\end{displaymath}
where $h,B,Q,R$ are measurable functions, $h,B$ being positive.
The statistic $S(x_1,\dots,x_m)=\sum R(x_i)$ is complete and sufficient.
\end{prop}
We apply the criterion with
$h(x)=\frac{x^{k-1}}{\Gamma(k)}$, $Q(\theta)=-\theta$, 
$R(\mathbf{x})=x$, $B(\theta)=1$, we obtain:
\begin{coro}
Under the asymptotic model, $S=\sum_i \Gamma^{(i)}$ is a sufficient and complete
statistic.
\end{coro}

\begin{prop}[Lehmann-Scheff\'e's Theorem]\label{Th:LS}[see \cite{Lind} Th.4, Chap.8]
Let $S$ be a complete and sufficient statistic, and let  $\xi^*$ be an unbiased
estimator.
The estimator $\mathbb{E}[\xi^*|S]$ is the unbiased estimator
with the lowest variance.
It is said to be \emph{efficient}.
\end{prop}
\begin{coro}
Let $\tilde{\xi}$ be an unbiased estimator of $\theta$.
\begin{equation}\label{Eq:LS} 
\mathbb{E}[(\tilde{\xi}-\theta)^2]\geq \mathbb{E}[(\hat{\xi}-\theta)^2]=\frac{\theta^2}{km-2}.
\end{equation}
\end{coro}
This is a consequence of Lehmann-Scheff\'e's Theorem with $S=\sum_{i=1}^m \Gamma^{(i)}$ and
$\xi^*=\hat{\xi}$, since
$
\mathbb{E}[\hat{\xi}|S]=\hat{\xi}.
$
This inequality is sharp for our model, but it is valid only for unbiased estimators.
Cramr-Rao inequality \cite[Th. 6.4, page 122]{Le} gives a more general lower bound :
\begin{prop}[Cramr-Rao inequality] 
Let $g_\theta$ be the density of $\Gamma^{(1)}$. Assume that $\theta\mapsto \log g_\theta(x)$ is
continuously differentiable, and that the quantity
\begin{displaymath}
I(\theta)=-\mathbb{E}\left[\frac{d^2}{d\theta^2}\log g_\theta(\Gamma^{(1)})\right]
\end{displaymath}
is finite and positive.
Let $\xi^*$ be a square-integrable function such that
$b(\theta):=\mathbb{E}[\xi^*]-\theta$ is continuously differentiable. then
\begin{displaymath}
\mathbb{E}\left[(\xi^*-\theta)^2\right]\geq \frac{(1+b'(\theta))^2}{mI(\theta)}
+b(\theta)^2.
\end{displaymath}
In particular, if $\xi^*$ is unbiased, its variance is bounded from below by $1/mI(\theta)$.
\end{prop}
The quantity $I(\theta)$ is the \emph{Fisher information}. In our case, it reduces to
\begin{align*}
I(\theta)&=-\mathbb{E}\left[\frac{d^2}{d\theta^2}\log
\left(\frac{(\Gamma^{(1)})^{k-1}\ \theta^{k}}{\Gamma(k)}\ e^{-\theta \Gamma^{(1)}}\right)\right],\\
&=-\mathbb{E}\left[\frac{d^2}{d\theta^2}\left(k\log\theta -\theta\Gamma^{(1)}\right)\right],\\
&=\frac{k}{\theta^2}.
\end{align*}
Under the asymptotic model, Cramr-Rao inequality reads
\begin{equation}\label{Eq:CR}
\mathbb{E}[(\xi^*-\theta)^2]
\geq (1+b'(\theta))^2\frac{\theta^2}{km}+b(\theta)^2,
\end{equation}
for any estimator $\xi^*$ such that $\theta\mapsto\mathbb{E}[\xi^*]$ is continuously differentiable. Cramr-Rao inequality confirms that quadratic error cannot decrease faster than $\theta^2/M$. The estimator $\hat{\xi}$ achieves this lower bound, up to a factor $km/(km-2)$.

\subsection{Proof of Proposition \ref{Prop:Conv}}\label{Sec:Preuve}
Let $A$ be the event
\begin{align*}
A=A_{k,m,\theta}\ =&\{\text{For any }1\leq i\leq m, \text{ at least $k$ hashed values
 lie in the $i$-th interval}\}\\
=&\{\text{For any }1\leq i\leq m, X_{(k),i}<\tfrac{1}{m}\}
\end{align*}
When $\theta\gg 2m^2$, $A$ occurs with a high probability :
\begin{align}
1-\mathbb{P}(A)&=\mathbb{P}\left(\cup_{i=1}^m \{\text{less than $k$ hashed values lie
in }[\tfrac{i-1}{m};\tfrac{i}{m})\}\right)\notag\\
&\leq m\ \mathbb{P}\left(\{\text{less than $k$ values lie in }
[0;\tfrac{1}{m})\}\right)\notag\\
&\leq m\ \mathbb{P}(\mathcal{B}_{\theta,1/m}<k)\notag\\
&\leq m\ \mathbb{P}(\mathcal{B}_{\theta,1/m}-\tfrac{\theta}{m}<-\tfrac{\theta}{2m})\notag\\
&\leq m\ \exp(-\theta/(2m^2)),\label{Eq:Hoeff}
\end{align}
in which $\mathcal{B}_{\theta,1/m}$ follows the binomial distribution with parameters $(\theta,1/m)$, 
and in which \eqref{Eq:Hoeff} follows from Hoeffding's inequality.
Now, write
\begin{displaymath}\label{Eq:EvA}
\mathbb{E}\left[\left(f(X_{(k),1},\ldots,X_{(k),m})-\theta\right)^2\right]
= \mathbb{E}[(f-\theta)^2\mathbf{1}_A]
+\mathbb{E}[(f-\theta)^2\mathbf{1}_{\overline{A}}].
\end{displaymath}
The restriction to $A$ of the distribution of $(X_{(k),1},\ldots,X_{(k),m})$
admits a density on  $\bbR^m$, that can be computed along the lines of  \cite[pages 8-13]{Dav}, leading to  :
\begin{multline*}
\mathbb{E}[\left(f(X_{(k),1},\ldots,X_{(k),m})-\theta\right)^2\mathbf{1}_A]=\\
\int_{[0,1/m]^m}(f(x)-\theta)^2\ \frac{\theta !\,(1-\textstyle{\sum_i x_i})^{\theta-mk}}{(\theta-km)!}\ \prod_{i=1}^m\frac{{x_i}^{k-1}}{\Gamma(k)}\ dx.
\end{multline*}
On the set $\overline{A}$, at least one of the $X_{(k),i}$'s is equal to
$1/m$, so that, using  \eqref{Eq:Hoeff} and the fact that $f$ is bounded outside a neighbourhood of $0$, we find
\begin{equation}\label{odetheta}
\mathbb{E}[(f-\theta)^2\mathbf{1}_{\overline{A}}]
\leq 2\left(\sup_{\norme{x}\geq 1/m} |f(x)|^2+\theta^2\right)(1-\mathbb{P}(A))\leq 4m\theta^2 \exp(-\theta/(2m)^2),
\end{equation}
which gives the last term in \eqref{Eq:Cv:Prop}, provided $\theta$ is large enough.
Thus  the proof reduces to show that for any $\eps>0$, 


$$
I=\mathcal{O}\left(\theta^{\eps-1}\mathbb{E}\left[\left(f(X_{(k),1},\ldots,X_{(k),m})-\theta\right)^2\right]\right),
$$
in which $I$ is defined below :
\begin{eqnarray*}
I
&=&
I_{1}-I_{2},
\\
I_{1}
&=&
\int_{[0,1/m]^m}\left(f(x)-\theta\right)^2\ \frac{\theta !\,(1-\textstyle{\sum_i x_i})^{\theta-mk}}{(\theta-km)!}\ \prod_{i=1}^m\frac{{x_i}^{k-1}}{\Gamma(k)}\ dx,
\\
I_{2}
&=&
\int_{\mathbb{R}_{+}^m}\left(f(x)-\theta\right)^2\ \theta^{km}\exp(-\theta \sum_i x_i)\ \prod_{i=1}^m \frac{{x_i}^{k-1}}{\Gamma(k)}\ dx.
\end{eqnarray*}
By the substitution $y_i=\theta x_i$, and with the notation $s=\sum_i y_i$, we get
\begin{align}
I_{1}
&=
\int_{[0,\theta/m]^m}\left(f\left(y/\theta\right)-\theta\right)^2\ \tfrac{\theta! \left(1-(s/\theta)\right)^{\theta-mk}}{(\theta-km)!\theta^{km}}\ \prod_{i\leq m} \frac{{y_i}^{k-1}}{\Gamma(k)}\ dy,
\label{Eq:HypNouv}
\\
I_{2}
&=
\int_{\mathbb{R}_{+}^m}\left(f\left(y/\theta\right)-\theta\right)^2\ \exp(-s)\ \prod_{i=1}^m \frac{{y_i}^{k-1}}{\Gamma(k)}\ dy.
\end{align}
Set
\begin{displaymath}
F(y,\theta)
\ =
\ 1\ -\ \frac{(\theta-km)!\ \theta^{km}\ \exp(-s)}{\theta! \left(1-\tfrac{s}{\theta}\right)^{\theta-mk}}
.
\end{displaymath}
We write $I=J-K$, with
\begin{align*}
J
&=
\int_{[0,\theta/m]^m}\ F(y,\theta)\ \left(f\left(y/\theta\right)-\theta\right)^2\ \tfrac{\theta! \left(1-(s/\theta)\right)^{\theta-mk}}{(\theta-km)!\theta^{km}}\ \prod_{i\leq m} \frac{{y_i}^{k-1}}{\Gamma(k)}\ dy,
\\
K
&=
\int_{\mathbb{R}_+^m\backslash \left[0,\theta/m\right]^m}\left(f\left(y/\theta\right)-\theta\right)^2\ \exp(-s)\ \prod_{i=1}^m \frac{{y_i}^{k-1}}{\Gamma(k)}\ dy.
\end{align*}
We will show that these two integrals are $\mathrm{o}(\theta^2)$. 
Since $f$ is continuous and bounded away from zero outside some neighbourhood of 0, there exist $c,c'>0$ such that in  $\mathbb{R}_+^m\backslash \left[0,\frac{\theta}{m}\right]^m$, when $\theta$ is large enough,
\begin{displaymath}
|f(y/\theta)-\theta|^2
\leq c+\theta^2\leq c'\theta^2\label{Eq:MajPoly}.
\end{displaymath}
Then
\begin{displaymath}
|K|\leq\ c'\theta^2\ \int_{\mathbb{R}_+^m\backslash \left[0,\frac{\theta}{m}\right]^m}\ \prod_{i\leq m} \frac{{y_i}^{k-1}}{\Gamma(k)}\ e^{-s}\ dy\,
\end{displaymath}
that vanishes exponentially fast, 
the integral on the right hand side being the probability that the maximum of an $m$-sample of Gamma distributed random variables is larger than $\theta/m$. Let us write $J=J_1+J_2$ in which
\begin{align*}
J_1
&=
\int_{[0,\theta^\alpha/m]^m}\ F(y,\theta)\ \left(f\left(y/\theta\right)-\theta\right)^2\ \tfrac{\theta! \left(1-(s/\theta)\right)^{\theta-mk}}{(\theta-km)!\theta^{km}}\ \prod_{i\leq m} \frac{{y_i}^{k-1}}{\Gamma(k)}\ dy,
\\
J_2&=\int_{\left[0,\frac{\theta}{m}\right]^m \backslash \left[0,\frac{\theta^\alpha}{m}\right]^m}\ F(y,\theta)\ \left(f\left(y/\theta\right)-\theta\right)^2\ \tfrac{\theta! \left(1-(s/\theta)\right)^{\theta-mk}}{(\theta-km)!\theta^{km}}\ \prod_{i\leq m} \frac{{y_i}^{k-1}}{\Gamma(k)}\ dy.
\end{align*}
and assume that $\alpha\in (0,1/2)$. Let
\[F(y,\theta)
=
1-e^{\psi(s,\theta)},
\]
in which
\begin{align*}
\psi(s,\theta)
&=
-s-\theta\ln\left(1-\tfrac{s}{\theta}\right)+mk\ln\left(1-\tfrac{s}{\theta}\right)-\sum_{\ell=1}^{km-1}\ln\left(1-\tfrac{\ell}{\theta}\right)
\\
&=
\mathcal{O}\left(\theta^{2\alpha-1}\right)+\mathcal{O}\left(\theta^{\alpha-1}\right)+\mathcal{O}\left(\theta^{-1}\right),
\end{align*}
as long as $y \in \left[0,\theta^\alpha/m\right]^m$. Using \eqref{Eq:HypNouv}, we conclude 
that there exists $c_\alpha >0$ such that for, $\theta$ large enough,
$$
J_1\leq c_\alpha \theta^{2\alpha -1}I_1 \leq c_\alpha \theta^{2\alpha -1}\mathbb{E}\left[\left(f(X_{(k),1},\ldots,X_{(k),m})-\theta\right)^2\right].
$$

As concerns $J_2$, we write $J_2=J_{2,1}-J_{2,2}$, with
\begin{align*}
J_{2,1}&=\int_{\left[0,\frac{\theta}{m}\right]^m \backslash \left[0,\theta^\alpha/m\right]^m}\ \left(f\left(y/\theta\right)-\theta\right)^2\ e^{-s}\ \prod_{i\leq m} \frac{{y_i}^{k-1}}{\Gamma(k)}dy,\\
J_{2,2}&=\int_{\left[0,\frac{\theta}{m}\right]^m \backslash \left[0,\theta^\alpha/m\right]^m}
\ \left(f\left(y/\theta\right)-\theta\right)^2\ \tfrac{\theta! \left(1-(s/\theta)\right)^{\theta-mk}}{(\theta-km)!\theta^{km}}\ \prod_{i\leq m} \frac{{y_i}^{k-1}}{\Gamma(k)}\ dy.
\end{align*}
First, we bound $J_{2,1}$ : by the second assumption of  Proposition \ref{Prop:Conv}, 
there exist $c,c'>0$ such that, for all
$y\in \left[0,\frac{\theta}{m}\right]^m \backslash \left[0,\tfrac{\theta^\alpha}{m}\right]^m$, 
\begin{align*}
|f\left(y/\theta\right)-\theta|^2
&\leq \theta^2+ f^2\left(y/\theta\right)\\
&\leq \theta^2 +c\ \frac{1}{1\wedge\norme{y/\theta}^{2r}}\\
&\leq \frac{c'\theta^{2r}}{\norme{y}^{2r}},
\end{align*}
in which $r$ can be chosen larger than 1. Since $y\in\left[0,\frac{\theta}{m}\right]^m \backslash \left[0,\tfrac{\theta^\alpha}{m}\right]^m$,
we have $s\geq \theta^\alpha/m$. It follows that
\begin{align*}
J_{2,1}
&\leq
\int_{\left[0,\frac{\theta}{m}\right]^m \backslash \left[0,\theta^\alpha/m\right]^m}\frac{c'\theta^{2r}\ e^{-s}}{\norme{y}^{2r}}\ \prod_{i\leq m} \frac{{y_i}^{k-1}}{\Gamma(k)}\ dy.
\\
&\leq
\ c'\theta^{2r}\ e^{-\theta^\alpha/m}\int_{\left[0,\frac{\theta}{m}\right]^m \backslash \left[0,\theta^\alpha/m\right]^m}\ \prod_{i\leq m} \frac{{y_i}^{k-1}}{\Gamma(k)}\ \frac{dy}{\norme{y}^{2r}},
\end{align*}
in which the last integral is polynomial in $\theta$. For $J_{2,2}$, we have, as soon as $\theta\ge 2mk$, 
\begin{align*}
\tfrac{\theta! \left(1-(s/\theta)\right)^{\theta-mk}}{(\theta-km)!\theta^{km}}\ 
&\leq
\ \left(1-(s/\theta)\right)^{\theta/2}
\ \leq\ \exp(-s/2).
\end{align*}
Thus, by the same argument,  $J_{2,2}$ vanishes exponentially fast.
To finish the proof, take $\eps =2\alpha$.

\section{Optimality of $\hat{\xi}$}\label{Sec:Reel}
We return to the original model, in which $X_{(k),i}$ denotes the $k$-th smallest value
lying in $[\frac{i-1}{m},\frac{i}{m})$.
We wish to discuss the optimality of
\begin{displaymath}
\hat{\xi}=\frac{km-1}{\sum_{i=1}^{m} X_{(k),i}}.
\end{displaymath}
The combination of \eqref{Eq:LS} and \eqref{Eq:CR} gives the following result, which
is the main result of the present work.
\begin{theo}[Optimality of $\hat{\xi}$]\label{Th:GrosTheo}
Let $\tilde{\xi}=\tilde{\xi}(\mathbf{x})$ be a continuous function on $\bbR_+^m-\set{0}$.
Assume that there exists $r>0$ such that, in the neighbourhood of $0$, 
$$
|f(x)| =\mathcal{O}\left(\frac{1}{\parallel x\parallel^r}\right).
$$
The application
$$
b(\theta):=\mathbb{E}[\tilde{\xi}(\Gamma^{(1)},\dots,\Gamma^{(m)})]-\theta,
$$
is continuously differentiable, and
$$
\mathbb{E}[(\tilde{\xi}(X_{(k),1},\dots,X_{(k),m})-\theta)^2]
\geq \frac{\theta^2}{km}\left(1+b'(\theta)\right)^2+\mathrm{o}(\theta^2).
$$
If we assume furthermore that $\tilde{\xi}$ is unbiased in the asymptotic model \ie 
\begin{equation}\label{Eq:HypSansBiais}
\mathbb{E}[\tilde{\xi}(\Gamma_1,\dots,\Gamma_m)]=\theta,
\end{equation}
then
\begin{align*}
\mathbb{E}[(\tilde{\xi}(X_{(k),1},\dots,X_{(k),m})-\theta)^2]
&\geq  \mathbb{E}[(\hat{\xi}-\theta)^2]+\mathrm{o}(\theta^2),\\
&=\frac{\theta^2}{km-2} +\mathrm{o}(\theta^2).
\end{align*}
\end{theo}
\begin{rem} The second assumption ($\tilde{\xi}$ unbiased in the asymptotic model) was implicitly made in \cite{Gi}.
\end{rem}
\begin{proof}
Proposition \ref{Prop:Conv} with $\eps=1/2$ implies that there exists $c>0$ such that, for $\theta$
large enough,
\begin{align*}
\mathbb{E}[(\tilde{\xi}(X_{(k),1},\dots,X_{(k),m})-\theta)^2]
&\geq\mathbb{E}[\tilde{\xi}(\Gamma_1,\dots,\Gamma_m)]
-c \theta^{-1/2}\mathbb{E}[(\tilde{\xi}(X_{(k),1},\dots,X_{(k),m})-\theta)^2]+\mathrm{o}(\theta^2),\\
&\geq \frac{\theta^2}{km}\left(1+b'(\theta)\right)^2
-c \theta^{-1/2}\mathbb{E}[(\tilde{\xi}(X_{(k),1},\dots,X_{(k),m})-\theta)^2]+\mathrm{o}(\theta^2)
\end{align*}
by \eqref{Eq:CR}.\LG{Petite modif de la preuve pour tre plus clair}
The conclusion of the Theorem follows : if $\mathbb{E}[(\tilde{\xi}-\theta)^2]$ is larger than, say, $\theta^{7/3}$ then
there is nothing to prove ; if it is smaller then the right-hand term is 
$$
\frac{\theta^2}{km}\left(1+b'(\theta)\right)^2+\mathrm{o}(\theta^2).
$$
If we assume furthermore that \eqref{Eq:HypSansBiais} holds,
then \eqref{Eq:LS} gives the second assertion of the theorem.
\end{proof}

\subsection{The case $m=1$}
The lower bound given by Theorem \ref{Th:GrosTheo} depends on the choice of
$(k,m)$ only through the product $km$. Thus, regardless of algorithmic
considerations, the quadratic error does not benefit from the partition of $[0,1]$ in $m$ sub-intervals, and we may assume $m=1$ to study the optimality of our estimator. 
When  $m=1$, the law of the observation $X_{(k),1}$ is easy to deal with,
and we obtain a sharp and exact lower bound (\ie  valid for any $\theta$).

The irrelevance, with respect to statistics,  of splitting $[0,1]$ into $m$ subintervals, is perhaps clearer in the asymptotic model. We noted in Section \ref{sec:asympt} that $\Xi_{k,m}$ is distributed, asymptotically, as $m$ independent copies of $T_{k}Y$, in which $Y$ is a $k$-sample of the exponential distribution with expectation $1/\theta$, and
\[T_{k}Y=(Y_{1},Y_{1}+Y_{2},Y_{1}+Y_{2}+Y_{3}, \dots, Y_{1}+Y_{2}+\dots+Y_{k}).\] 
From a statistical perspective, since $T_{k}$ is one-to-one, there is no loss of information from $Y$ to $T_{k}Y$. Thus, for the estimation of $\theta$, $\Xi_{k,m}$ is equivalent to $m$ independent copies of a $k$-sample of the exponential distribution with expectation $1/\theta$, i.e.  a $km$-sample of the exponential distribution with expectation $1/\theta$. But this is also equivalent to   $T_{km}\hat{Y}$, in which $\hat{Y}$ is a $km$-sample of the exponential distribution with expectation $1/\theta$ :   we can obtain $T_{km}\hat{Y}$ as the asymptotic distribution of the first $km$ order statistics of a $\theta$-sample of uniform random variables on $[0,1]$, that is, without splitting the interval $[0,1]$ into $m$ subintervals.  For the exponential distribution, the sum of the $km$ elements of the sample is known to be  complete and sufficient : when splitting, this corresponds to the sum of the $m$ copies of the $k$-th order statistic, and when not splitting, this sum is asymptotic to the $km$-th order statistic. 
\begin{theo}
Consider algorithm \textsc{MinCount}, with parameter $m=1$. It returns
\begin{displaymath}
\hat{\xi}(X_{(k)})=\frac{k-1}{X_{(k)}}.
\end{displaymath}
For any $\theta$, $\hat{\xi}$ is an unbiased estimator. Furthermore, $\hat{\xi}$ is the unbiased
estimator with the lowest variance.
\end{theo}
Note that in this particular case $m=1$, our estimator coincides with the estimator $\xi_{1}$ proposed by Giroire.
\LG{Petite remarque ajout\'ee}
\begin{proof}
We know from \eqref{Eq:EstNaif} that $\hat{\xi}$ is unbiased.
We apply the  Lehmann-Scheff\'e Theorem to the law $Q_\theta$ of $X_{(k)}$.
\begin{align*}
Q_\theta (x)dx&=\theta\binom{\theta-1}{k-1}x^{k-1}(1-x)^{\theta-k}dx\\
&=B(\theta)h(x)\exp((\theta -k)\log (1-x))dx,
\end{align*}
with notations of Proposition \ref{Prop:Critere}.
This shows that the statistic $\log (1-X_{(k)})$ is complete and sufficient. Thus the variance  of the statistic
\begin{displaymath}
\mathbb{E}\left[\left.\frac{k-1}{X_{(k)}}\ \right|\ \log (1-X_{(k)})\right]\ =\ \frac{k-1}{X_{(k)}}
\end{displaymath}
is minimal.
\end{proof}

\bibliographystyle{plain}
\bibliography{BiblioEstimationNickel}
\end{document}